\documentstyle[amscd,amssymb,verbatim,epsf]{amsart}

\begin{document}
\theoremstyle{plain}
\newtheorem{Thm}{Theorem}
\newtheorem{Cor}{Corollary}
\newtheorem{Con}{Conjecture}
\newtheorem{Main}{Main Theorem}
\newtheorem{Lem}{Lemma}
\newtheorem{Prop}{Proposition}

\theoremstyle{definition}
\newtheorem{Def}{Definition}
\newtheorem{Note}{Note}

\theoremstyle{remark}
\newtheorem{notation}{Notation}
\renewcommand{\thenotation}{}

\errorcontextlines=0
\numberwithin{equation}{section}
\renewcommand{\rm}{\normalshape}%

\title[On maximal surfaces ]%
   {On maximal surfaces in the space of oriented geodesics of hyperbolic 3-space}
\author{Nikos Georgiou}
\address{Nikos Georgiou\\
          Department of Computing and Mathematics\\
          Institute of Technology, Tralee\\
          Clash\\
          Tralee\\
          Co. Kerry\\
          Ireland.}
\email{nikos.georgiou@@research.ittralee.ie}

\keywords{Kaehler structure, hyperbolic 3-space, area stationary surfaces}
\subjclass{Primary: 51M09; Secondary: 51M30}

\date{10th February 2010}

\begin{abstract}
We study area-stationary, or maximal, surfaces in the space ${\mathbb L}({\mathbb H}^3)$ of oriented geodesics of hyperbolic 3-space, 
endowed with the canonical neutral K\"ahler structure. We prove that every holomorphic curve in ${\mathbb L}({\mathbb H}^3)$ is a maximal surface. 
We then classify Lagrangian maximal surfaces $\Sigma$ in ${\mathbb L}({\mathbb H}^3)$ and prove that the family of parallel surfaces in 
${\mathbb H}^3$ orthogonal to the geodesics $\gamma\in \Sigma$ form a family of equidistant tubes around a geodesic.
\end{abstract}

\maketitle

The last two decades has seen increasing interest in spaces ${\mathbb L}({\mathbb M})$ of oriented geodesics of a manifold ${\mathbb M}$, with
particular attention to their rich geometric structure. In the case of the space ${\mathbb L}({\mathbb E}^3)$ of 
oriented affine lines in Euclidean 3-space 
this interest can be traced back over a hundred years to Weierstrass's construction of minimal surfaces \cite{weierstrass} and  Whittaker's 
solutions to the wave equation \cite{whittaker}. 

A natural complex structure on ${\mathbb L}({\mathbb E}^3)$ was considered by Hitchin to 
construct monopoles in ${\mathbb E}^3$ \cite{hitch}, and then Guilfoyle and Klingenberg understood that the canonical symplectic structure 
on ${\mathbb L}({\mathbb E}^3)$ is compatible with this complex structure \cite{klbg}\cite{klbg1} and that the associated K\"ahler metric is of neutral
signature. Salvai subsequently proved that this neutral K\"ahler metric is (up to addition of the round metric)
the unique metric on ${\mathbb L}({\mathbb E}^3)$ that is invariant under Euclidean motions \cite{salvai0}. This K\"ahler structure has recently 
been used  by Guilfoyle and Klingenberg to solve an 80 year old conjecture of Carath\'eodory \cite{caratheodory}.  

More recently, Anciaux, Guilfoyle and Romon have studied maximal Lagrangian surfaces in $TN$, with $N$ being an oriented Riemannian surface and the 
neutral K\"ahler structure generalising that of the space of oriented geodesics in Euclidean and Lorentzian 3-space.

In addition, Salvai constructed a neutral K\"ahler metric on the space ${\mathbb L}({\mathbb H}^3)$ of oriented geodesics in hyperbolic 3-space 
\cite{salvai}, while the geometry of ${\mathbb L}({\mathbb H}^3)$ was explored by Georgiou, Guilfoyle and Klingenberg \cite{gag} \cite{gagak} 
\cite{nikosbrendan}.

Recently Alekseevsky, Guilfoyle and Klingenberg have given a description of all metrics defined on the space ${\mathbb L}({\mathbb M})$ of oriented 
geodesics of space forms or rank 1 riemannian symmetric spaces, that are invariant under the isometry group of ${\mathbb M}$ 
\cite{Alekseevsky}.

The current paper can be considered as a continuation of the work of Georgiou and Guilfoyle (see \cite{nikosbrendan}) on surface theory 
of ${\mathbb L}({\mathbb H}^3)$. Here we study immersed maximal surfaces, that is, immersed surfaces that are stationary for the area functional. 
These surfaces have locally maximum area with respect to area variations, because of the neutral signature $(++--)$ of the metric defined on the 
ambient space ${\mathbb L}({\mathbb H}^3)$. 

Since ${\mathbb L}({\mathbb H}^3)=S^2\times S^2-\{diag\}$, we can consider surfaces in  ${\mathbb L}({\mathbb H}^3)$ of rank 0, 1 or 2, depending on
the rank of the projection onto the first factor.

In the first two sections we illustrate the geometric background on the construction of ${\mathbb L}({\mathbb H}^3)$ and its
submanifold theory, including geodesics and surfaces (further details can be found in \cite{gag}\cite{gagak} \cite{nikosbrendan}).

In section 3 we investigate rank 1 surfaces in ${\mathbb L}({\mathbb H}^3)$ and prove that there are no holomorphic curves of rank 1. Moreover, in the
case of a surface being Lagrangian of rank 1, we prove that it can not be maximal.

Surfaces of rank 2 are studied in section 4 and, in particular, we prove that every holomorphic curve of rank 2 is maximal and, since the only holomorphic 
curve of rank 0 are orthogonal to a horosphere, we conclude that every holomorphic curve is maximal.  

We also classify all Lagrangian maximal surfaces. We summarize this result as:

\vspace{0.1in}

\noindent {\bf Main Theorem}:

{\it Let $S\subset{\mathbb{H}}^3$ be a $C^3$ smooth immersed oriented surface and $\Sigma\subset{\mathbb{L}}({\mathbb{H}}^3)$ 
be the Lagrangian surface formed by the oriented geodesics normal to $S$. 

The surface $\Sigma$ is maximal iff $S$ is an equidistant tube around a geodesic in ${\mathbb{H}}^3$. In terms of holomorphic coordinates 
$(\mu_1,\mu_2)$ on ${\mathbb L}({\mathbb H}^3)$, the surface $\Sigma$ is given by
\[
\bar{\mu}_2=\frac{1+\lambda_1\mu_1}{\lambda_2+\mu_1}
\]
where $\lambda_1,\lambda_2\in{\mathbb C}$ with  $\lambda_1\lambda_2\neq 1$.
}

\vspace{0.2in}

\section{The Neutral K\"ahler Metric on ${\mathbb{L}}({\mathbb{H}}^3)$}

We briefly recall the basic construction of the canonical neutral K\"ahler metric on the space 
${\mathbb{L}}({\mathbb{H}}^3)$ of oriented geodesics of ${\mathbb{H}}^3$ - further details can be found in \cite{gag}.
We use one of two models of ${\mathbb H}^3$, the Poincar\'e ball model:
\[
 B^3=\{y\in {\mathbb{R}}^3\;|\; |y|^2< 1\},
\]
and hyperbolic metric
\[
d\tilde{s}^2=\frac{4}{(1-|y|^2)^2}|dy|^2,
\]
where $|\cdot|$ is the Euclidean norm, while for the upper-half space model:
\[
{\mathbb{R}}^3_{+}=\{\;(x_{0},x_{1},x_{2})\in {\mathbb{R}}^3 \;|\; x_{0}>0\;\},
\]
with standard coordinates $(x_{0},x_{1},x_{2})$ on ${\mathbb{R}}^3$, the hyperbolic metric 
has expression:
\begin{equation}\label{e:hypmet}
ds^2=\frac{(dx_{0})^2+(dx_{1})^2+(dx_{2})^2}{(x_{0})^2}.
\end{equation}
These are related by the mapping 
${\mathbb{R}}^3_{+}\rightarrow  B^3\colon (x_{0},x_{1},x_{2}) \mapsto (y_{1},y_{2},y_{3})$
defined by
\[
y_{1}=\frac{2x_1}{(x_0+1)^2+(x_1)^2+(x_2)^2},\quad 
y_{2}=\frac{2x_2}{(x_0+1)^2+(x_1)^2+(x_2)^2},
\]
\[
y_{3}=\frac{(x_0)^2+(x_1)^2+(x_2)^2-1}{(x_0+1)^2+(x_1)^2+(x_2)^2}.
\]

An oriented geodesic in ${\mathbb{H}}^3$ is uniquely determined by its beginning and end point on the boundary of the 
ball model, and so ${\mathbb{L}}({\mathbb{H}}^3)$ can be identified with $S^2\times S^2-\Delta$, where $\Delta$ is the 
diagonal in $S^2\times S^2$. Endowing $S^2\times S^2-\Delta$ with the standard differentiable structure, a tangent 
vector to an oriented geodesic $\gamma\in {\mathbb{L}}({\mathbb{H}}^3)$ can then be identified with an orthogonal 
Jacobi field along $\gamma\subset {\mathbb{H}}^3$. 

Rotation of Jacobi fields through 90$^{0}$ about $\gamma$ defines an almost complex structure on 
${\mathbb{L}}({\mathbb{H}}^3)$. This almost complex structure is integrable, and so ${\mathbb{L}}({\mathbb{H}}^3)$ 
becomes a complex surface, which turns out to be biholomorphic to 
${\mathbb{P}}^1\times {\mathbb{P}}^1-\overline{\Delta}$. Here $\overline{\Delta}$ is the ``reflected'' diagonal: 
in terms of holomorphic coordinates $(\mu_1,\mu_2)$ on 
${\mathbb{P}}^1\times {\mathbb{P}}^1$, $\overline{\Delta}=\{(\mu_1,\mu_2):\mu_1\bar{\mu}_2=-1\}$. 

The complex structure ${\mathbb{J}}$ on ${\mathbb{L}}({\mathbb{H}}^3)$ can be supplemented with a compatible symplectic 
structure $\Omega$, which has the following expression in holomorphic coordinates:
\begin{equation}\label{e:sympl}
\Omega=-\left[\frac{1}{(1+\mu_{1}\bar{\mu}_{2})^2}d\mu_{1}\wedge d\bar{\mu}_{2}
           +\frac{1}{(1+\bar{\mu}_{1}\mu_{2})^2}d\bar{\mu}_{1}\wedge d\mu_{2}\right].
\end{equation}
Together we obtain a K\"ahler metric ${\mathbb{G}}(\cdot\; ,\cdot)=\Omega({\mathbb{J}}\cdot\; ,\cdot)\;\;$:
\begin{equation}\label{e:metric}
{\mathbb{G}}=-i\left[\frac{1}{(1+\mu_{1}\bar{\mu}_{2})^2}d\mu_{1}\otimes d\bar{\mu}_{2}
                   -\frac{1}{(1+\bar{\mu}_{1}\mu_{2})^2}d\bar{\mu}_{1}\otimes d\mu_{2}\right].
\end{equation}
This metric, which has signature $++--$, is invariant under the action induced on ${\mathbb{L}}({\mathbb{H}}^3)$ by 
the isometry group of ${\mathbb{H}}^3$. Indeed, this has been shown to be the unique K\"ahler metric on 
${\mathbb{L}}({\mathbb{H}}^3)$ with this property \cite{salvai}.  

In order to transfer geometric data between ${\mathbb{L}}({\mathbb{H}}^3)$ and ${\mathbb{H}}^3$ we use a correspondence 
space:

\begin{figure} [ht]
\vspace*{0.9cm}

\unitlength0.5cm

\begin{picture}(16,6)
\put(5.1,4.8){$\pi_1$}
\put(5.4,6.5){${\mathbb{L}}({\mathbb{H}}^3)\times{\mathbb{R}}$}
\put(7,6){\vector(1,-1){2.8}}
\put(8.5,4.8){$\Phi$}
\put(5.2,2){${\mathbb{L}}({\mathbb{H}}^3)$}
\put(6,6){\vector(0,-1){3}}
\put(10,2.4){${\mathbb{H}}^3$}
\end{picture}
\end{figure}

The key property of this correspondence is that, given $\gamma\in{\mathbb L}({\mathbb H}^3)$, the set 
$\Phi\circ\pi_1^{-1}(\gamma)$ is the oriented geodesic in ${\mathbb H}^3$, while, for a point $p\in{\mathbb H}^3$,
$\pi_1\circ\Phi^{-1}(p)$ is the set of oriented geodesics in ${\mathbb L}({\mathbb H}^3)$ that pass through $p$.

The map $\Phi$ takes an oriented geodesic $\gamma$ in ${\mathbb L}({\mathbb H}^3)$ and a real number $r$ to
the point on  $\gamma$ an affine parameter distance $r$ from some fixed point on the geodesic. This choice of
point on each geodesic can be made globally, but we more often just use a local choice, which is sufficient for our
purposes.

In terms of holomorphic coordinates ($\mu_1,\mu_2$) on ${\mathbb L}({\mathbb H}^3)$ and upper-half space
coordinates ($x_0,x_1,x_2$) the map $\Phi$ has expression:
\begin{equation}\label{e:phi}
z=\frac{1-\mu_1\bar{\mu}_2}{2\bar{\mu}_2}+\left(\frac{1+\mu_1\bar{\mu}_2}{2\bar{\mu}_2}\right)\tanh r,
\qquad\qquad 
t=\frac{|1+\bar{\mu}_1\mu_2|}{2|\mu_2|\cosh r}.
\end{equation}
where $z=x_1+ix_2$ and $t=x_0$.

\vspace{0.2in}

\section{Surfaces in ${\mathbb{L}}({\mathbb{H}}^3)$}

Our interest in this paper is focused on the study of two parameter families of oriented geodesics, or surfaces in ${\mathbb L}({\mathbb H}^3)$. 
Therefore, we recall some basic results on the surface theory of the space of oriented geodesics in hyperbolic 3-space - further details can be found 
in \cite{nikosbrendan}.

For computational purposes, we give explicit local parameterizations of the geodesic congruence. The dual picture of a surface in 
${\mathbb L}({\mathbb H}^3)$ is to consider the surfaces in ${\mathbb H}^3$ that are orthogonal to a given set of geodesics. However, not every geodesic 
congruence has such orthogonal surfaces - indeed, most don't. To explain this further, we consider the first order properties of $\Sigma$, 
which can be described by two complex functions, the optical scalars: $\rho,\sigma:\Sigma\times {\mathbb R}\mapsto {\mathbb C}$. 
The real part $\Theta$ and imaginary part $\lambda$ of $\rho$ are the {\it divergence} and {\it twist} of the geodesic congruence, 
while $\sigma$ is the {\it shear}.

\begin{Def}
A {\it null frame} in ${\mathbb{H}}^3$ is a trio $\{e_{0},e_{+},e_{-}\}$ of complex vector fields in 
${\mathbb{C}}\otimes $T${\mathbb{H}}^3$ where $e_{0}$ is real, $e_{+}$ is the complex conjugate of $e_{-}$ 
and they satisfy the following properties:
\[
<e_{0}\; ,\; e_{0}>=1,\quad <e_{0}\; ,\; e_{+}>=<e_{+}\; ,\; e_{+}>=0,\quad <e_{+}\; ,\; e_{-}>=1,
\]
where $<\; ,\; >$ is the hyperbolic metric.
\end{Def}
Given an orthonormal frame $\{e_0,e_1,e_2\}$ on T${\mathbb{H}}^3$ we construct a null frame by
\[
e_{+}=\frac{1}{\sqrt{2}}(e_1-ie_2),\qquad e_{-}=\frac{1}{\sqrt{2}}(e_1+ie_2).
\]
\begin{Def}
Given a surface $\Sigma\subset {\mathbb{L}}({\mathbb{H}}^3)$ an {\it adapted null frame} is a null frame 
$\{e_{0},e_{+},e_{-}\}$ such that, for each $\gamma\in\Sigma$, we have $e_{0}=\dot{\gamma}$, and the orientation 
of $\{e_{0},{\mathbb{R}}{\mbox{e}}(e_{+}),{\mathbb{I}}{\mbox{m}}(e_{+})\}$ is the standard orientation on 
${\mathbb{H}}^3$.

For a given surface $\Sigma$ in ${\mathbb{L}}({\mathbb{H}}^3)$ and a null frame, the optical scalars can be defined by:
\[
\rho=<\nabla_{e_{0}} e_{+}, e_{-}>\qquad\sigma=<\nabla_{e_{0}} e_{+}, e_{+}>,
\]
where $<,>$ is the hyperbolic metric with hyperbolic connection $\nabla$.
\end{Def}
In terms of the holomorphic coordinates $(\mu_1,\mu_2)$, the optical scalars have the following local expressions:
\begin{equation}\label{e:0sigma0}
\sigma=\frac{8\mu_2 J_{\bar{2}\;\bar{1}}}{\bar{\mu}_2\Delta |1+\mu_1\bar{\mu}_2|^2},
\end{equation}
\begin{equation}\label{e:0rho0}
\rho=-1-\frac{8e^{-r}}{\Delta}\left[\frac{J_{2\bar{1}}}{(1+\bar{\mu}_1\mu_2)^2}e^{r}-\frac{|\mu_2|^2 J_{1\bar{1}}}{|1+\bar{\mu}_1\mu_2|^2}e^{-r}\right],
\end{equation}
where
\begin{equation}\label{e:jacobians11}
J_{kl}=\partial\mu_{k}\bar{\partial}\mu_{l}-\bar{\partial}\mu_{k}\partial\mu_{l}\qquad k,l=1,2,\bar{1},\bar{2},
\end{equation}
and 
\begin{equation}\label{e:1delta1}
\frac{1}{4}\Delta=\frac{J_{2\bar{2}}}{|\mu_2|^2|1+\mu_1\bar{\mu}_2|^2}e^{2r}+\frac{J_{\bar{2}1}}{(1+\mu_1\bar{\mu}_2)^2}+\frac{J_{\bar{1}2}}{(1+\bar{\mu}_1\mu_2)^2}+\frac{|\mu_2|^2 J_{1\bar{1}}}{|1+\mu_1\bar{\mu}_2|^2}e^{-2r}.
\end{equation}

\begin{Def}
A surface $\Sigma$ in ${\mathbb L}({\mathbb H}^3)$, is said to be {\it Lagrangian} if the symplectic form $\Omega$ pulled back to $\Sigma$ vanishes. 
\end{Def} 
In order to avoid any confusion, whether we use the expression of {\it geodesic congruence} we mean a surface in ${\mathbb L}({\mathbb H}^3)$. 

We now give the following important Proposition:
\begin{Prop}
 \cite{nikosbrendan} The following statements are equivalent:
\begin{enumerate}
\item[(i)] the geodesic congruence $\Sigma$ is Lagrangian,
\item[(ii)] locally there exists a surface $S$ in ${\mathbb H}^3$ such that the geodesics of $\Sigma$ are normal to $S$,
\item[(iii)] the imaginary part $\lambda$ of $\rho$ (the twist) is zero.
\end{enumerate}
\end{Prop}
Let $\Sigma$ be a Lagrangian surface in ${\mathbb L}({\mathbb H}^3)$ parameterized by $\nu\mapsto (\mu_1(\nu,\bar{\nu}),\mu_2(\nu,\bar{\nu}))$. 
The surfaces $S$ in ${\mathbb H}^3$ orthogonal to the geodesics of  $\Sigma$ are given by equations (\ref{e:phi}), where the functions 
$r = r(\nu,\bar{\nu})$ solve:
\begin{equation}\label{e:itgivesthelagrangiancond}
2\partial r=\frac{\mu_2}{\bar{\mu}_1\mu_2+1}\left(\partial\bar{\mu}_1+\frac{\partial\mu_2}{\mu^2_2}\right)+\frac{\bar{\mu}_2}{\mu_1\bar{\mu}_2+1}\left(\partial\mu_1+\frac{\partial\bar{\mu}_2}{\bar{\mu}^2_2}\right),
\end{equation}
where $\partial$ denotes the derivative with respect of $\nu$.

The set of Lagrangian geodesic congruences is divided into three categories, depending on the rank of the immersion of the geodesic congruence.
\begin{Def}
Given an immersion $f:\Sigma\rightarrow {\mathbb L}({\mathbb H}^3)$, consider the map $(\pi\circ f)_{\ast}:T\Sigma\rightarrow T{\mathbb P}^1$, where $\pi$ is projection onto the first factor of ${\mathbb L}({\mathbb H}^3)={\mathbb{P}}^1\times {\mathbb{P}}^1-\overline{\Delta}$. The {\it rank} of the immersion $f$ at a point $\gamma\in\Sigma$ is defined to be the rank of this map at $\gamma$, which can be 0, 1 or 2.
\end{Def}

Note that by reversing the orientation of the geodesics, the rank can be defined
by projection onto the second factor.
A rank 0 Lagrangian geodesic congruence correspond to a geodesic congruence orthogonal to a horosphere.

In the Lagrangian case, the functions $\sigma$ and $\rho$ have the following interpretation in terms of the
second fundamental form of the orthogonal surfaces in ${\mathbb H}^3$.

\begin{Prop}\label{p:lag}
\cite{nikosbrendan} Let $S\subset {\mathbb{H}}^3$ be a $C^2$ immersed surface and $\Sigma\subset {\mathbb{L}}({\mathbb{H}}^3)$ 
be the oriented normal geodesic congruence. Then
\begin{equation}\label{e:optical}
|\sigma|=\frac{1}{2}|\lambda_1-\lambda_2|\qquad\qquad \rho=-\frac{1}{2}(\lambda_1+\lambda_2),
\end{equation}
where $\lambda_1$ and $\lambda_2$ are the principal curvatures of $S$.
\end{Prop}

The induced metric ${\mathbb G}_{\Sigma}$ on a Lagrangian surface $\Sigma$ in ${\mathbb{L}}({\mathbb{H}}^3)$ can be described by the functions $\sigma$ and $\lambda$:
\begin{Thm}
\cite{nikosbrendan} Let $\Sigma$ be a surface in ${\mathbb{L}}({\mathbb{H}}^3)$. The induced metric is Lorentz (degenerate, 
Riemannian) iff $\; |\sigma|^2-\lambda^2>0\; (=0,<0)$, where $\lambda$ and $\sigma$ are the twist and the shear of $\Sigma$.
\end{Thm}

The theorem says that if a surface is Lagrangian then is either Lorentz or degenerate, the latter occurring at umbilic points on the orthogonal 
surfaces in ${\mathbb H}^3$.

The following theorem recovers the Weierstrass representation for the flat surfaces in ${\mathbb{H}}^3$:
\begin{Thm} \cite{nikosbrendan} 
Let $S$ be an oriented $C^2$ smooth immersed surface in ${\mathbb{H}}^3$ with normal geodesic congruence $\Sigma$. Assume that the metric 
${\mathbb G}_{\Sigma}$ induced on $\Sigma$ by the neutral K\"ahler metric is non-degenerate.

Then S is flat iff $\Sigma$ is of rank two and is parameterized by $\mu_2=\mu_2(\bar{\mu}_1)$, that is, $\mu_2$ is an anti-holomorphic function of $\mu_1$.
\end{Thm}

We use the complex structure ${\mathbb J}$ of ${\mathbb{L}}({\mathbb{H}}^3)$ in order to describe another important class of surfaces:

\begin{Def}
A point $\gamma$ on a surface $\Sigma\subset{\mathbb{L}}({\mathbb{H}}^3)$ is said to be a {\it complex point} 
if the complex structure ${\mathbb{J}}$ acting on ${\mathbb{L}}({\mathbb{H}}^3)$ preserves $T_{\gamma}\Sigma$. A surface
$\Sigma\subset{\mathbb{L}}({\mathbb{H}}^3)$ is said to be a {\it holomorphic curve} if all of the points of $\Sigma$ are
complex points.
\end{Def}
In particular:
\begin{Prop}\label{p:sigmavanishing}
\cite{nikosbrendan} A point $\gamma$ on a surface $\Sigma$ is complex iff the shear vanishes along $\gamma$.
\end{Prop}
By Proposition \ref{p:lag}, observe that complex points on a Lagrangian surface $\Sigma\subset{\mathbb{L}}({\mathbb{H}}^3)$ correspond to 
umbilic points on the surfaces in ${\mathbb{H}}^3$ orthogonal to $\Sigma$.

\vspace{0.2in}

\section{Non-existence of rank one maximal surfaces}

Consider a surface $\Sigma\subset {\mathbb L}({\mathbb H}^3)$ of rank one. Then $\Sigma$ can be locally parameterized by $\mu_1=\mu_1(s)$ and 
$\mu_2=\mu_2(s,t)$ where $(s,t)\in D$ with $D$ being an open subset of ${\mathbb R}^2$. 

We prove first that $\Sigma$ cannot be a holomorphic curve:

\begin{Prop}\label{p:nocomplexpoints1}
An immersed holomorphic curve in ${\mathbb L}({\mathbb H}^3)$ cannot be of rank 1.
\end{Prop}
\begin{pf}
Assume the existence of an immersed rank one surface on ${\mathbb L}({\mathbb H}^3)$ such that in an open neighborhood $U\subset \Sigma$ is 
holomorphic, 
which, by Proposition \ref{p:sigmavanishing}, is equivalent to the vanishing of the shear $\sigma$ in this open set. Considering now the 
local parametrisation $\Sigma\rightarrow {\mathbb L}({\mathbb H}^3):(s,t)\mapsto (\mu_1(s),\mu_2(s,t))$, the vanishing of the shear implies the 
vanishing of $J_{\bar{2}\;\bar{1}}$, which means $\partial_{s}\bar{\mu}_1\partial_{t}\bar{\mu}_2=0$ on $U$. 

Then, in an open subset $V$ of $U$, either $\partial_{s}\bar{\mu}_1=0$ or $\partial_{t}\bar{\mu}_2=0$. In any case, we have a contradiction since, 
if $\partial_{s}\bar{\mu}_1=0$ then $\mu_1$ is constant and therefore $\Sigma$ is of rank 0 on $V$, and if $\partial_{t}\bar{\mu}_2=0$ then 
$\mu_2$ only depends on $s$ and $\Sigma$ would not be immersed.
\end{pf}

We now assume that the rank 1 surface $\Sigma$ is Lagrangian. In this case the induced metric $g=f^{\ast}{\mathbb G}$ has components in local coordinates $(s,t)$:
\[
g_{ss}=2{\mbox Im}\left[\frac{\partial_{s}\mu_1\partial_{s}\bar{\mu}_2}{(1+\mu_1\bar{\mu}_2)^2}\right]\qquad
g_{st}={\mbox Im}\left[\frac{\partial_{s}\mu_1\partial_{t}\bar{\mu}_2}{(1+\mu_1\bar{\mu}_2)^2}\right]\qquad
g_{tt}=0,
\]
and the nonzero Christoffel symbols are given by:
\[
\Gamma^{s}_{ss}={\mathbb R}{\mbox{e}}\left(\frac{\partial_{s}^2\mu_1}{\partial_{s}\mu_1}-\frac{2\bar{\mu}_2\partial_{s}\mu_1}{1+\mu_1\bar{\mu}_2}\right)\qquad 
\Gamma^{t}_{st}={\mathbb R}{\mbox{e}}\left(\frac{\partial_{st}^2\mu_2}{\partial_{t}\mu_2}-\frac{2\bar{\mu}_1\partial_{s}\mu_2}{1+\bar{\mu}_1\mu_2}\right),
\]
\[
\Gamma^{t}_{tt}={\mathbb R}{\mbox{e}}\left(\frac{\partial_{t}^2\mu_2}{\partial_{t}\mu_2}-\frac{2\bar{\mu}_1\partial_{t}\mu_2}{1+\bar{\mu}_1\mu_2}\right).
\]
It is already known that the induced metric $g$ of a rank one Lagrangian surface $\Sigma$ is scalar flat \cite{nikosbrendan}.

The second fundamental form $h=h_{ij}^{\mu_k}$ has non-vanishing components:  
\[
h_{ss}^{\mu_1}=\partial_{s}^2\mu_1-\frac{2\bar{\mu}_2(\partial_{s}\mu_1)^2}{1+\mu_1\bar{\mu}_2}-\partial_{s}\mu_1\Gamma_{ss}^{s},
\quad 
h_{ss}^{\mu_2}=\partial_{s}^2\mu_2-\frac{2\bar{\mu}_1(\partial_{s}\mu_2)^2}{1+\bar{\mu}_1\mu_2}-\partial_{s}\mu_2\Gamma_{ss}^{s}-\partial_{t}\mu_2\Gamma_{ss}^{t}
\]
\[
h_{st}^{\mu_2}=\partial_{st}^2\mu_2-\frac{2\bar{\mu}_1\partial_{s}\mu_2\partial_{t}\mu_2}{1+\bar{\mu}_1\mu_2}-\partial_{t}\mu_2\Gamma_{st}^{t},\qquad h_{tt}^{\mu_2}=\partial_{t}^2\mu_2-\frac{2\bar{\mu}_1(\partial_{t}\mu_2)^2}{1+\bar{\mu}_1\mu_2}-\partial_{t}\mu_2\Gamma_{tt}^{t},
\]
with $h_{ij}^{\bar{\mu}_{k}}=\overline{h_{ij}^{\mu_{k}}}$.

\vspace{0.1in}

Unlike with the space ${\mathbb L}({\mathbb E}^3)$ of oriented lines in Euclidean 3-space \cite{wilhelbrendanareastationary}, the following 
Proposition shows that there are no maximal Lagrangian surface in ${\mathbb L}({\mathbb H}^3)$ of rank one:
\begin{Prop}\label{p:nolagrangian1}
There are no maximal Lagrangian surfaces in ${\mathbb L}({\mathbb H}^3)$ of rank one. 
\end{Prop}
\begin{pf}
Let $\Sigma$ be a Lagrangian surface in ${\mathbb L}({\mathbb H}^3)$ of rank one, locally parameterized by $\mu_1=\mu_1(s)$ and $\mu_2=\mu_2(s,t)$.

Firstly, we find the mean curvature vector $H=2{\mathbb R}e(H^{\mu_1}\partial/\partial\mu_1+H^{\mu_2}\partial/\partial\mu_2)$ in local coordinates $(s,t)$. The components $H^{\mu_{i}}$ are given by 
\[
H^{\mu_{i}}=g^{ss}h^{\mu_{i}}_{ss}+2g^{st}h^{\mu_{i}}_{st}+g^{tt}h^{\mu_{i}}_{tt}
\]
Then $H^{\mu_1}=0$ and it remains to find $H^{\mu_2}$. By using the expressions of $h^{\mu_{i}}_{ij}$ and by considering the Lagrangian condition:
\begin{equation}\label{e:lagrangianrank1}
\frac{\partial_{s}\mu_1\partial_{t}\bar{\mu}_2}{(1+\mu_1\bar{\mu}_2)^2}=-\frac{\partial_{s}\bar{\mu}_1\partial_{t}\mu_2}{(1+\bar{\mu}_1\mu_2)^2},
\end{equation}
we find $H^{\mu_2}$, and finally the mean curvature vector $H$ of $\Sigma$ is:
\[
H=4{\mathbb R}e\left[g^{st}\frac{(1+\bar{\mu}_1\mu_2)^2}{\partial_{s}\bar{\mu}_1}\partial_{t}\left(\frac{\partial_{s}\mu_1\partial_{s}\bar{\mu}_2}{(1+\mu_1\bar{\mu}_2)^2}+\frac{\partial_{s}\bar{\mu}_1\partial_{s}\mu_2}{(1+\bar{\mu}_1\mu_2)^2}\right)\frac{\partial}{\partial\mu_2}\right],
\]
which means that the surface $\Sigma$ is area stationary iff
\[
\partial_{t}\left[{\mathbb R}e\left(\frac{\partial_{s}\mu_1\partial_{s}\bar{\mu}_2}{(1+\mu_1\bar{\mu}_2)^2}\right)\right]=0.
\]
The above condition and the Lagrangian condition give:
\begin{equation}\label{e:maxcondit2}
A\partial_{s}\bar{\mu}_2+\bar{A}\partial_{s}\mu_2=f(s), 
\end{equation}
\begin{equation}\label{e:maxcondit3}
A\partial_{t}\bar{\mu}_2+\bar{A}\partial_{t}\mu_2=0,
\end{equation}
where
\[
A=\frac{\partial_{s}\mu_1}{(1+\mu_1\bar{\mu}_2)^2}.
\]
Differentiate equations (\ref{e:maxcondit2}) and (\ref{e:maxcondit3}) with respect to $t$ and $s$, respectively, and then subtract:
\begin{equation}\label{e:condjhgf}
{\mathbb R}e(\partial_{t}A\partial_{s}\bar{\mu}_2-\partial_{s}\bar{A}\partial_{t}\bar{\mu}_2)=0.
\end{equation}
After a brief computation we get:
\begin{align}
\partial_{s}A\partial_{t}\bar{\mu}_2&=\frac{\partial_{s}^2\mu_1\partial_{t}\bar{\mu}_2}{(1+\mu_1\bar{\mu}_2)^2}-\frac{2\mu_1\partial_{s}\mu_1\partial_{s}\bar{\mu}_2\partial_{t}\bar{\mu}_2}{(1+\mu_1\bar{\mu}_2)^3}-\frac{2\bar{\mu}_2(\partial_{s}\mu_1)^2\partial_{t}\bar{\mu}_2}{(1+\mu_1\bar{\mu}_2)^3},\nonumber \\
\partial_{t}A\partial_{s}\bar{\mu}_2&=-\frac{2\mu_1\partial_{s}\mu_1\partial_{s}\bar{\mu}_2\partial_{t}\bar{\mu}_2}{(1+\mu_1\bar{\mu}_2)^3},\nonumber
\end{align}
and then condition (\ref{e:condjhgf}) becomes
\begin{equation}\label{e:maxcondit3jk}
{\mathbb R}e\left(\frac{2\bar{\mu}_2(\partial_{s}\mu_1)^2\partial_{t}\bar{\mu}_2}{(1+\mu_1\bar{\mu}_2)^3}-\frac{\partial_{s}^2\mu_1\partial_{t}\bar{\mu}_2}{(1+\mu_1\bar{\mu}_2)^2}\right)=0.
\end{equation}
Using the Lagrangian condition (\ref{e:lagrangianrank1}) in (\ref{e:maxcondit3jk}), we have
\begin{equation}\label{e:maxcondit3}
\frac{\bar{\mu}_2\partial_{s}\mu_1}{1+\mu_1\bar{\mu}_2}-\frac{\mu_2\partial_{s}\bar{\mu}_1}{1+\bar{\mu}_1\mu_2}=\frac{1}{2}\left(\frac{\partial_{s}^2\mu_1}{\partial_{s}\mu_1}-\frac{\partial_{s}^2\bar{\mu}_1}{\partial_{s}\bar{\mu}_1}\right)=h(s).
\end{equation}
Integration of (\ref{e:lagrangianrank1}) with respect of $t$ gives 
\begin{equation}\label{e:maxcondit4}
\frac{\bar{\mu}_2\partial_{s}\mu_1}{1+\mu_1\bar{\mu}_2}+\frac{\mu_2\partial_{s}\bar{\mu}_1}{1+\bar{\mu}_1\mu_2}=g(s),
\end{equation}
and then the sum (\ref{e:maxcondit3})+(\ref{e:maxcondit4}), is
\[
\frac{\bar{\mu}_2\partial_{s}\mu_1}{1+\mu_1\bar{\mu}_2}=h(s)+g(s)=m(s)\partial_{s}\mu_1.
\]
Hence
\[
\bar{\mu}_2=\frac{m}{1-m\mu_1}=\bar{\mu}_2(s),
\]
which is a contradiction, since $\Sigma$ is of rank one. Therefore there are no Lagrangian maximal surfaces of rank one.
\end{pf}

\vspace{0.2in}

\section{Rank two maximal surfaces}

Consider a rank 2 surface $\Sigma$ in ${\mathbb L}({\mathbb H}^3)$. That is, a surface $\Sigma$ given locally by $\mu_1\rightarrow (\mu_1,\mu_2(\mu_1,\bar{\mu}_1))$ for some smooth function $\mu_2:{\mathbb C}\rightarrow {\mathbb C}$. We are interested in maximal surfaces in ${\mathbb L}({\mathbb H}^3)$ of rank 2 and therefore we consider variations of the area integral
\[
{\cal A}(\Sigma)=\int_\Sigma |{\mathbb G}|^{\scriptstyle{\frac{1}{2}}}d\mu_1 d\bar{\mu}_1.
\]
For an arbitrary parameterization $\mu_1\rightarrow (\mu_1,\mu_2(\mu_1, \bar{\mu}_1))$ the area integral is
\[
|{\mathbb G}|=\frac{\Delta^2}{64}(\lambda^2-|\sigma|^2),
\]
where
\[
\lambda=\frac{4i}{\Delta}\left[\frac{\partial\mu_2}{(1+\bar{\mu}_1\mu_2)^2}-\frac{\bar{\partial}\bar{\mu}_2}{(1+\mu_1\bar{\mu}_2)^2}\right],\quad
\sigma=\frac{8\mu_2 \partial\bar{\mu}_2}{\bar{\mu}_2\Delta |1+\mu_1\bar{\mu}_2|^2},
\]
\[
\frac{1}{4}\Delta=\frac{\partial\mu_2\bar{\partial}\bar{\mu}_2-\bar{\partial}\mu_2\partial\bar{\mu}_2}{|\mu_2|^2|1+\mu_1\bar{\mu}_2|^2}e^{2r}-\frac{\bar{\partial}\bar{\mu}_2}{(1+\mu_1\bar{\mu}_2)^2}-\frac{\partial\mu_2}{(1+\bar{\mu}_1\mu_2)^2}+\frac{|\mu_2|^2}{|1+\mu_1\bar{\mu}_2|^2}e^{-2r},
\]
with $\partial$ denotes the differentiation with respect to $\mu_1$.

\vspace{0.1in}

A surface is maximal if $\delta {\cal A}(\Sigma)=0$. In order to compute this quantity note that
\[
\frac{\Delta^2\lambda^2}{64}=-\frac{1}{4}\left[\frac{\partial \mu_2}{(1+\bar{\mu}_1\mu_2)^2}-\frac{\bar{\partial}\bar{\mu}_2}{(1+\mu_1\bar{\mu}_2)^2}\right]^2,
\]
and so
\[
\delta\left(\frac{\Delta^2\lambda^2}{64}\right)=-{\mbox Re}\left[\frac{\partial \mu_2}{(1+\bar{\mu}_1\mu_2)^2}-\frac{\bar{\partial}\bar{\mu}_2}{(1+\mu_1\bar{\mu}_2)^2}\right]
     \left[\frac{\partial \delta \mu_2}{(1+\bar{\mu}_1\mu_2)^2}-\frac{2\bar{\mu}_1\partial \mu_2\delta \mu_2}{(1+\bar{\mu}_1\mu_2)^3}\right],
\]
while, since
\[
\frac{\Delta^2|\sigma|^2}{64}=\frac{\partial\bar{\mu}_2\bar{\partial}\mu_2}{|1+\bar{\mu}_1\mu_2|^4},
\] 
we have 
\[
\delta\left(\frac{\Delta^2|\sigma|^2}{64}\right)=2{\mbox Re}\left(\frac{\partial\bar{\mu}_2\bar{\partial}\delta \mu_2}{|1+\bar{\mu}_1\mu_2|^4}
              -\frac{2\bar{\mu}_1}{1+\bar{\mu}_1\mu_2}\frac{\partial\bar{\mu}_2\bar{\partial}\mu_2\delta \mu_2}{|1+\bar{\mu}_1\mu_2|^4}\right).
\]
Combining these we find that
\begin{align}
\delta|{\mathbb G}|^{\scriptstyle{\frac{1}{2}}}&=\frac{16}{\Delta\sqrt{\lambda^2-|\sigma|^2}}{\mbox Re}\left\{-\frac{1}{2}\left[\frac{\partial \mu_2}{(1+\bar{\mu}_1\mu_2)^2}       -\frac{\bar{\partial}\bar{\mu}_2}{(1+\mu_1\bar{\mu}_2)^2}\right]\partial \left(\frac{\delta \mu_2}{(1+\bar{\mu}_1\mu_2)^2}\right)    \right.\nonumber\\
&\qquad\qquad\qquad\qquad\qquad\qquad\qquad-\left.\frac{\partial\bar{\mu}_2\bar{\partial}\delta \mu_2}{|1+\bar{\mu}_1\mu_2|^4}
              +\frac{2\bar{\mu}_1}{1+\bar{\mu}_1\mu_2}\frac{\partial\bar{\mu}_2\bar{\partial}\mu_2\delta \mu_2}{|1+\bar{\mu}_1\mu_2|^4}\right\}\nonumber\\
=&2{\mbox Re}\left[\frac{\lambda i}{\sqrt{\lambda^2-|\sigma|^2}}\partial \left(\frac{\delta \mu_2}{(1+\bar{\mu}_1\mu_2)^2}\right)
    -\frac{\bar{\mu}_2\sigma\bar{\partial}\delta \mu_2}{\mu_2|1+\bar{\mu}_1\mu_2|^2\sqrt{\lambda^2-|\sigma|^2}}\right.\nonumber\\
&\qquad\qquad\qquad\qquad\qquad\qquad\qquad +\left.\frac{16\bar{\mu}_1}{1+\bar{\mu}_1\mu_2}\frac{\partial\bar{\mu}_2\bar{\partial}\mu_2\delta\mu_2}{\Delta|1+\bar{\mu}_1\mu_2|^4\sqrt{\lambda^2-|\sigma|^2}}\right].\nonumber
\end{align}
Integrating by parts we have established the following:
\begin{Prop}\label{p:maximall01}
A rank two surface is maximal iff
\begin{align}
\frac{-i}{(1+\bar{\mu}_1\mu_2)^2}\partial \left(\frac{\lambda}{\sqrt{\lambda^2-|\sigma|^2}}\right)&
    +\bar{\partial}\left(\frac{\bar{\mu}_2\sigma}{\mu_2|1+\bar{\mu}_1\mu_2|^2\sqrt{\lambda^2-|\sigma|^2}}\right)+\nonumber\\
&\qquad\qquad\qquad\qquad    +\frac{\bar{\mu}_1\Delta|\sigma|^2}{4(1+\bar{\mu}_1\mu_2)\sqrt{\lambda^2-|\sigma|^2}}=0.\nonumber
\end{align}
\end{Prop}

The following Proposition shows that all holomorphic curves on ${\mathbb L}({\mathbb H}^3)$ are maximal:

\begin{Prop}\label{p:maximall02}
Every holomorphic curve $\Sigma$ in ${\mathbb L}({\mathbb H}^3)$, where the metric ${\mathbb G}_{\Sigma}$ induced on $\Sigma$ by the neutral K\"ahler metric being non-degenerate, is maximal.
\end{Prop}
\begin{pf}
Consider a holomorphic curve $\Sigma$ in ${\mathbb L}({\mathbb H}^3)$. Then by Proposition \ref{p:sigmavanishing} the shear $\sigma$ vanishes throughout the surface $\Sigma$.

By Proposition \ref{p:nocomplexpoints1} we know that a holomorphic surface $\Sigma$ can be either rank 0 or 2. In the case of rank 0, the surface $\Sigma$ is totally null and, in particular, it is orthogonal to a horosphere, which is not our case. Then $\Sigma$ must be of rank 2 and therefore Proposition \ref{p:maximall01} shows that is maximal.
\end{pf}

Consider now a Lagrangian surface $\Sigma$ in ${\mathbb L}({\mathbb H}^3)$. We are interested in maximal Lagrangian surfaces of rank 2. In this case, the twist $\lambda$ vanishes on $\Sigma$ and then Proposition \ref{p:maximall01} implies that a Lagrangian surface $\Sigma$ of rank 2 will be maximal iff
\[
\partial\ln\left(\frac{\bar{\sigma}_0}{\sigma_0}\right)-\frac{4\bar{\mu}_2}{1+\mu_1\bar{\mu}_2}=0,
\]
where 
\begin{equation}\label{e:tosigma00}
\sigma_0=\frac{\partial\bar{\mu}_2}{(1+\mu_1\bar{\mu}_2)^2}\; .
\end{equation}
\begin{Def}
The {\it Lagrangian angle} $\phi$ of the surface $\Sigma\subset{\mathbb L}({\mathbb H}^3)$ is defined by
\[
\sigma_0=|\sigma_0|e^{2i\phi},
\]
where $\sigma_0$ is given by (\ref{e:tosigma00}).
\end{Def}

An equivalent condition that characterizes Lagrangian maximal surfaces in ${\mathbb L}({\mathbb H}^3)$ is given by the following Proposition:

\begin{Prop}\label{p:importantproposition1}
Let $\Sigma\subset {\mathbb L}({\mathbb H}^3)$ be a Lagrangian surface of rank two. Then $\Sigma$ is a maximal surface iff $\Sigma$ is locally orthogonal to a flat surface in ${\mathbb H}^3$ and the Lagrangian angle $\phi$, satisfies the following PDE:
\begin{equation}\label{e:finalequationforminimal0}
e^{-i\phi}\partial^2 e^{-i\phi}= e^{i\phi}\bar{\partial}^2 e^{i\phi}=|\sigma_0|.
\end{equation}
\end{Prop}
\begin{pf}
Assume that $\Sigma$ is a Lagrangian area stationary surface of rank two. Then $H^{\mu_{1}}=0$ which means that 
\[
\partial\ln\left(\frac{\bar{\sigma}_0}{\sigma_0}\right)-\frac{4\bar{\mu}_2}{1+\mu_1\bar{\mu}_2}=0
\]
and by introducing the Lagrangian angle $\phi$, the above gives
\begin{equation}\label{e:sxesiphi2}
\mu_2=\frac{i\bar{\partial}\phi}{1-i\bar{\mu}_1\bar{\partial}\phi}.
\end{equation}
By derivation of the above with respect of $\mu_1$ we obtain
\[
i\partial\bar{\partial}\phi=\rho_0\qquad -i\partial\bar{\partial}\phi=\bar{\rho}_0.
\]
The Lagrangian condition $\rho_0=\bar{\rho}_0$ implies that $\rho_0=0$ and therefore $\mu_2$ is anti-holomorphic, which means that $\Sigma$ is locally orthogonal to a flat surface in ${\mathbb H}^3$.

\noindent Because of $\rho_0=0$ we obtain 
\begin{equation}\label{e:expressionofmuwithres1}
\partial\bar{\partial}\phi=0.
\end{equation}

The fact that $\mu_2$ is an anti-holomorphic function of $\mu_1$ implies that $\ln\bar{\sigma}_0$ is anti-holomorphic too, which means that $\partial\ln\bar{\sigma}_0=0$ and then
\begin{equation}\label{e:sxesiphi1}
\partial\ln|\sigma_0|=2i\partial\phi.
\end{equation}
The expression of $\sigma_0$ in terms of $\phi$ is 
\begin{equation}\label{e:sxesiphi3}
\sigma_0=\frac{\partial\bar{\mu}_2}{(1+\mu_1\bar{\mu}_2)^2}=-(\partial\phi)^2-i\partial^2\phi.
\end{equation}
Then we have
\[
|\sigma_0|e^{2i\phi}=-(\partial\phi)^2-i\partial^2\phi
\]
which gives
\[
|\sigma_0|e^{i\phi}=[-(\partial\phi)^2-i\partial^2\phi\; ]\; e^{-i\phi}
\]
and therefore $|\sigma_0|=e^{-i\phi}\partial^2 e^{-i\phi}$, which implies equation (\ref{e:finalequationforminimal0}).

\vspace{0.1in}

We now prove the converse. If $\phi$ is a real solution of (\ref{e:finalequationforminimal0}), it satisfies 
\begin{equation}\label{e:sxesiphi3000}
-i\partial^2\phi-(\partial\phi)^2=|\sigma_0|e^{2i\phi}=\sigma_0.
\end{equation}
By the assumption that $\Sigma$ is locally orthogonal to a flat in ${\mathbb H}^3$, $\mu_2$ is anti-holomorphic and therefore $\sigma_0$ is a 
holomorphic function. Then 
\[
\bar{\partial}[-i\partial^2\phi-(\partial\phi)^2]=0,
\]
which implies 
\[
\partial[(\bar{\partial}\partial\phi) e^{-2i\phi}]=0,
\]
and hence there is a holomorphic function $\beta$ such that 
\[
(\bar{\partial}\partial\phi) e^{-2i\phi}=\beta
\qquad
(\bar{\partial}\partial\phi) e^{2i\phi}=\bar{\beta}.
\]
Therefore, the function $\phi$ can be written as
\begin{equation}\label{e:sxesiphi300001}
\phi=a+\bar{a},
\end{equation}
where $a$ is a holomorphic function. In other words, we have proved that $\bar{\partial}\partial\phi=0$.

On the other hand, using equations (\ref{e:sxesiphi3000}) and (\ref{e:sxesiphi300001}), $\mu_2$ must satisfies the following equation:
\begin{equation}\label{e:sxesiphi3001}
-i\partial^2a-(\partial a)^2=\frac{\partial\bar{\mu}_2}{(1+\mu_1\bar{\mu}_2)^2},
\end{equation}
but because of the fact that $\bar{\mu}_2$ and $a$ are holomorphic, the above equation is equivalent to an ordinary differential equation of first order. 
In addition, we observe that
\[
\bar{\mu}_2=-\frac{i\partial a}{1+i\mu_1\partial a}
\]
is a solution of  (\ref{e:sxesiphi3001}) and because this equation is equivalent to an ODE of first order, it is unique.

Then it is easy to see that $H^{\mu_1}=0$ and therefore $\Sigma$ is a Lagrangian maximal surface.
\end{pf}
In the following Proposition we give an explicit local expression of all Lagrangian maximal surfaces in ${\mathbb L}({\mathbb H}^3)$ in terms of the 
holomorphic coordinates $(\mu_1,\mu_2)$ on ${\mathbb P}^1\times {\mathbb P}^1-\bar{\Delta}$:
\begin{Prop}
Every Lagrangian maximal surface $\Sigma$ in ${\mathbb L}({\mathbb H}^3)$ of rank two can be locally parameterized by 
\begin{equation}\label{e:maximalimmersion}
\Sigma\rightarrow {\mathbb L}({\mathbb H}^3):(\mu_1,\bar{\mu}_1)\mapsto \left(\mu_1,\mu_2=\frac{\bar{\lambda}_1\bar{\mu}_1+1}{\bar{\mu}_1+\bar{\lambda}_2}\right),
\end{equation}
where $\lambda_1,\lambda_2\in {\mathbb C}$ with  $\lambda_1\lambda_2\neq 1$.
\end{Prop}
\begin{pf}
Let $\Sigma$ be a Lagrangian maximal surface of rank two in ${\mathbb L}({\mathbb H}^3)$. By Proposition \ref{p:importantproposition1} 
the surface $\Sigma$ is locally orthogonal to a flat surface in ${\mathbb H}^3$, which allows us to obtain the holomorphic parameterization 
$(\mu_1,\bar{\mu}_1)\mapsto (\mu_1,\mu_2(\bar{\mu}_1))$. In addition, the Lagrangian angle $\phi$ must satisfies 
equation (\ref{e:finalequationforminimal0}).

There is a holomorphic function $a$ such that $\phi=a+\bar{a}$, and applying this to equation (\ref{e:finalequationforminimal0}), we get:
\[
e^{3i\bar{a}}\bar{\partial}^2 e^{i\bar{a}}=e^{-3ia}\partial^2 e^{-ia}=c_{0},
\]
where $c_0\in {\mathbb R}$ is a real constant.

Then the holomorphic function $a$ satisfies
\begin{align}\label{e:diaforikiexisosimea1}
\partial^2 e^{-ia}=c_0 e^{3ia},
\end{align}
which is equivalent to the following ordinary differential equation of second order:
\[
\ddot{x}=c_0x^{-3}.
\] 
The unique solution of (\ref{e:diaforikiexisosimea1}) is
\begin{align}\label{e:diaforikiexisosimea2}
a=\frac{i}{2}\log[(\alpha_0\mu_1+\beta_0)^2-c_0]-\frac{i}{2}\log\alpha_0,
\end{align}
and the Lagrangian angle is $\phi=a+\bar{a}$.

The immersion of $\Sigma$ is obtained  by substituting (\ref{e:diaforikiexisosimea2}) into (\ref{e:sxesiphi2}) and then
\[
\mu_2=\frac{i\bar{\partial}\phi}{1-i\bar{\mu}_1\bar{\partial}\phi}
=\frac{i\bar{\partial}\bar{a}}{1-i\bar{\mu}_1\bar{\partial}\bar{a}}=\frac{\alpha_0^2\bar{\mu}_1+\alpha_0\beta_0}
     {\alpha_0\beta_0\bar{\mu}_1+\beta_0^2-c_0}.
\]
If we set $\bar{\lambda}_1=\alpha_0\beta_0^{-1}$ and $\bar{\lambda}_2=(\beta_0^2-c_0)(\alpha_0\beta_0)^{-1}$ then the maximal surface $\Sigma$ is 
given by the immersion (\ref{e:maximalimmersion}).

If $\lambda_1\lambda_2=1$ we find that $\Sigma$ is a totally null surface given by the immersion $\mu_2=\lambda_1$, and so it is not of rank two. 
\end{pf}

For a given Lagrangian maximal surface $\Sigma$ in ${\mathbb L}({\mathbb H}^3)$, there is locally a family of parallel flat surfaces in ${\mathbb H}^3$ 
such that their oriented normals are contained in $\Sigma$. We recall the classification of complete flat surfaces in hyperbolic 3-space:

\begin{Prop}\label{p:flatsurf0}
\cite{sasaki} \cite{Volkov} Let $S$ be a complete flat surface in hyperbolic 3-space ${\mathbb H}^3$. Then $S$ is either a horosphere or an equidistant 
tube of a geodesic in ${\mathbb H}^3$.
\end{Prop}

To proof the main theorem we need to introduce a particular class of surface in  hyperbolic 3-space ${\mathbb H}^3$:

\begin{Def}
A surface $S$ in hyperbolic 3-space ${\mathbb H}^3$ is called {\it isoparametric} if the principal curvatures of $S$ are constant.
\end{Def}

Note that all parallel surfaces $\{S_{t}\}_{t\in I}$ to the isoparametric surface $S$ are also isoparametric. 

The following Proposition gives a classification of the isoparametric surfaces in hyperbolic 3-space:

\begin{Prop}\label{p:isoparametric0}
\cite{Cartan} Let $S$ be an isoparametric surface in ${\mathbb H}^3$. Then $S$ is either a totally geodesic hyperbolic 2-space, or a totally umbilical surface or an equidistant tube around a geodesic.
\end{Prop}

We now prove our main result:

\vspace{0.1in}

\noindent {\bf Main Theorem}:

{\it Let $S\subset{\mathbb{H}}^3$ be a $C^3$ smooth immersed oriented surface and $\Sigma\subset{\mathbb{L}}({\mathbb{H}}^3)$ 
be the Lagrangian 
surface formed by the oriented geodesics normal to $S$.

The surface $\Sigma$ is maximal iff $S$ is an equidistant tube around a geodesic.}

\begin{pf}
Let $\Sigma$ be a Lagrangian geodesic congruence formed by the oriented geodesics normal to $S$. 

First assume that $\Sigma$ is maximal. Since it cannot be of rank 0 - as that would mean that it is totally null - and by Proposition 
\ref{p:nolagrangian1} it cannot be of rank 1, we conclude that $\Sigma$  is of rank two. Thus it is given locally by the graph: 
\[
\bar{\mu}_2=\frac{1+\lambda_1\mu_1}{\mu_1+\lambda_2},
\]
where $\lambda_1,\lambda_2\in {\mathbb C}$. 

The non-degeneracy condition of the induced  metric ${\mathbb G}_{\Sigma}$ implies that $\lambda_1\lambda_2\neq 1$.

In this case, an orthogonal surface $S\subset {\mathbb H}^3$ can be obtained by solving the following differential equation
\[
2\partial r=\frac{\partial\mu_2}{\mu_2(1+\bar{\mu}_1\mu_2)}+\frac{\partial\bar{\mu}_2}{\bar{\mu}_2(1+\mu_1\bar{\mu}_2)}+\frac{\bar{\mu}_2}{1+\mu_1\bar{\mu}_2},
\]
and, by using the fact that $\mu_2$ is holomorphic, we obtain, after a brief computation, that
\[
2\partial r=\frac{\lambda_1}{1+\lambda_1\mu_1},
\]
which implies
\begin{equation}\label{e:theformular0}
r=\frac{1}{2}\log|1+\lambda_1\mu_1|^2+r_0.
\end{equation}

The function $\Delta$, given by (\ref{e:1delta1}), is:
\[
\Delta=\frac{4|\lambda_1|^2[e^{-2r_0}-|\lambda_1\lambda_2-1|^2e^{2r_0}]}{|(\lambda_1\mu_1+1)^2+\lambda_1\lambda_2-1|^2}.
\]
The optical scalars $\rho$ and $\sigma$ of the Lagrangian maximal surface $\Sigma$ given by the expressions (\ref{e:0sigma0}) and (\ref{e:0rho0}) are:
\[
\sigma=\frac{2(\lambda_1\lambda_2-1)}{e^{-2r_0}-|\lambda_1\lambda_2-1|^2e^{2r_0}}\cdot\frac{1+\bar{\lambda}_1\bar{\mu}_1}{1+\lambda_1\mu_1}\qquad \mbox{and}\qquad\rho=-1+\frac{2}{1-e^{4r_0}|\lambda_1\lambda_2-1|^2}.
\]
If we denote by $h$ the mean curvature of the surface $S\subset {\mathbb H}^3$, Proposition \ref{e:optical} gives:
\begin{align}\label{e:constantmeancurvature0}
h=1+\frac{2}{e^{4r_0}|\lambda_1\lambda_2-1|^2-1}.
\end{align}

Consider now the principal curvatures  $m_1$ and $m_2$ of the surface $S$. The fact that $S$ is flat means that $m_1 m_2=1$. Then the mean curvature of the surface $S$ is 
\[
h=\frac{m_1+m_2}{2}=\frac{m_1+m_1^{-1}}{2},
\]
and by using the relation (\ref{e:constantmeancurvature0}), we observe that $m_1$ must satisfies the following quadratic equation 
\begin{align}\label{e:constantmeancurvature1}
m_1^2-2\left(1+\frac{2}{e^{4r_0}|\lambda_1\lambda_2-1|^2-1}\right)m_1+1=0.
\end{align}
Therefore the principal curvatures of the surface $S$ are constant and in particular are given by:
\[
m_1=\frac{e^{2r_0}|\lambda_1\lambda_2-1|+1}{e^{2r_0}|\lambda_1\lambda_2-1|-1},\qquad m_2=\frac{e^{2r_0}|\lambda_1\lambda_2-1|-1}{e^{2r_0}|\lambda_1\lambda_2-1|+1},
\]
and hence the surface $S$ is isoparametric. Propositions \ref{p:flatsurf0} and \ref{p:isoparametric0} tell us that the surface $S$ 
can be either a horosphere or an equidistant tube around a geodesic. By previous work (see the papers \cite{gagak} and \cite{nikosbrendan}) 
we have seen that geodesic congruences orthogonal to horospheres are totally null (the induced metric is degenerate). Therefore the surface 
$S$ must be an equidistant tube around a geodesic $\gamma$. 

In fact, every maximal surface $\Sigma$ is orthogonal to the set $\{S_{r_0}\}_{r_o\in {\mathbb R}}$ of all parallel equidistant tubes around a geodesic $\gamma$ and each such a surface $S_{r_0}$ is of hyperbolic distance $r_0$ from the surface $S_0$.

\vspace{0.1in}

Conversely, assume that the surface $S\subset {\mathbb H}^3$ is an equidistant tube around a geodesic 
$\gamma'$ with holomorphic coordinates $(\mu_1=\mu'_1,\mu_2=\mu'_2)$. Then 
$S$ belongs to the set of all parallel equidistant tubes ${\cal U}_{\gamma'}=\{S_{r_0}\}_{r_0\in {\mathbb R}}$ around the geodesic $\gamma'$. 
We first find an explicit expression of the orthogonal geodesic congruence $\Sigma\subset {\mathbb L}({\mathbb H}^3)$ to all surfaces in 
${\cal U}_{\gamma'}$.

Consider hyperbolic 3-space ${\mathbb{H}}^3$ in the half space model with coordinates $(x_0,x_1,x_2)$ and metric given by equation (\ref{e:hypmet}).

For a given point $p=(p_0,p_1,p_2)$ in ${\mathbb{H}}^3$ and a given vector 
$e_0=a_0\frac{\partial}{\partial x_0}+a_1\frac{\partial}{\partial x_1}+a_2\frac{\partial}{\partial x_2}\in T_{p}{\mathbb{H}}^3$ we now find the 
unique geodesic $\gamma:I\subset{\mathbb{R}}\rightarrow {\mathbb{H}}^3:r\mapsto \gamma(r)\in {\mathbb{H}}^3$ such that
\[
\gamma(0)=p\qquad\qquad \dot{\gamma}(0)=e_0,
\]
where $I$ is an open interval containing $0$ and the dot denotes the differentiation with respect of $r$.

Denote the Levi-Civita connection of $({\mathbb{H}}^3,ds^2)$ by $\nabla$. Then $\nabla_{\dot{\gamma}}\dot{\gamma}=0$ yields
\[
\ddot{x}_0+\frac{1}{x_0}(-\dot{x}_0^2+\dot{x}_1^2+\dot{x}_2^2)=0\qquad\qquad \ddot{x}_1-\frac{2}{x_0}\dot{x}_0\dot{x}_1=0\qquad\qquad \ddot{x}_2-\frac{2}{x_0}\dot{x}_0\dot{x}_2=0. 
\]
The first integrals are
\[
\frac{\dot{x}_0^2+\dot{x}_1^2+\dot{x}_2^2}{x_0^2}=\frac{a_0^2+a_1^2+a_2^2}{p_0^2} \qquad\qquad c_1=\frac{\dot{x}_1}{x_0^2}=\frac{a_1^2}{p_0^2}\qquad\qquad c_2=\frac{\dot{x}_2}{x_0^2}=\frac{a_2^2}{p_0^2}
\]
The geodesic $\gamma$ is
\[
x_0=p_0\sqrt{\frac{a_0^2+a_1^2+a_2^2}{a_1^2+a_2^2}}\cosh^{-1}\left[\frac{\sqrt{a_0^2+a_1^2+a_2^2}}{p_0}(r+r_0)\right]
\]
\[
x_1=\frac{a_1p_0\sqrt{a_0^2+a_1^2+a_2^2}}{a_1^2+a_2^2}\tanh\left[\frac{\sqrt{a_0^2+a_1^2+a_2^2}}{p_0}(r+r_0)\right]+c_3
\]
\[
x_2=\frac{a_2p_0\sqrt{a_0^2+a_1^2+a_2^2}}{a_1^2+a_2^2}\tanh\left[\frac{\sqrt{a_0^2+a_1^2+a_2^2}}{p_0}(r+r_0)\right]+c_4.
\]
The initial conditions $x_0(0)=p_0$ and $\dot{x}_0(0)=a_0$ yield
\[
\cosh\left(\frac{\sqrt{a_0^2+a_1^2+a_2^2}}{p_0}r_0\right)=\sqrt{\frac{a_0^2+a_1^2+a_2^2}{a_1^2+a_2^2}}
\]
\[
\sinh\left(\frac{\sqrt{a_0^2+a_1^2+a_2^2}}{p_0}r_0\right)=-\frac{a_0}{\sqrt{a_1^2+a_2^2}}
\]
thus,
\[
\tanh\left(\frac{\sqrt{a_0^2+a_1^2+a_2^2}}{p_0}r_0\right)=-\frac{a_0}{\sqrt{a_0^2+a_1^2+a_2^2}}
\]
Introduce complex coordinate $z=x_1+ix_2$ and set $t=x_0$. We then obtain
\[
\xi=c_1+ic_2=\frac{\beta}{t_0^2}\qquad\qquad \eta=c_3+ic_4=z_0+t_0\frac{a}{\bar{\beta}}, 
\]
where $t_0=t(0),\; z_0=z(0),\;\beta=a_1+ia_2$ and $a=a_0$.

Therefore for a given point $p=(z_0,t_0)$ and a given vector $e_0=a\frac{\partial}{\partial t}+\beta\frac{\partial}{\partial z}+\bar{\beta}\frac{\partial}{\partial\bar{z}}$ the unique oriented geodesic $\gamma=(\xi,\eta)$ with the initial conditions $\gamma(0)=p$ and $\dot{\gamma}(0)=e_0$ is given by
\begin{equation}\label{e:geodesic}
\xi=\frac{\beta}{t_0^2}\qquad\qquad \eta=z_0+t_0\frac{a}{\bar{\beta}}
\end{equation}

Fix the point $p$ on the given oriented geodesic $\gamma'=(\xi',\eta')$. Let $\gamma=(\xi,\eta)$ be an oriented geodesic that intersects 
$\gamma'$ orthogonally at $p$. Denote the unit tangent vectors of $\gamma,\gamma'$ at $p$ by $e_0,e'_0$ respectively. 
The orthogonality condition gives the following relation:
\[
e_0=\frac{1}{\sqrt{2}}(e^{-i\theta}e'_{+}+e^{i\theta}e'_{-}),
\] 
for some $\theta\in [0,2\pi )$ where
\[
e'_{+}=\frac{1}{\sqrt{2}|\xi'|\cosh^2 r_0}\frac{\partial}{\partial t}+\frac{1}{\sqrt{2}\cosh^2 r_0}\left(-\frac{e^{-r_0}}{\bar{\xi'}}\frac{\partial}{\partial z}+\frac{e^{r_0}}{\xi'}\frac{\partial}{\partial\bar{z}}\right)\qquad e'_{-}=\bar{e}'_{+}.
\]
Thus the unit tangent vector of $\gamma$ is
\[
e_0=\frac{\cos\theta}{|\xi'|\cosh^2 r_0}\frac{\partial}{\partial t}+\frac{\sinh(r_0+i\theta)}{\bar{\xi'}\cosh^2 r_0}\frac{\partial}{\partial z}+\frac{\sinh(r_0-i\theta)}{\xi'\cosh^2 r_0}\frac{\partial}{\partial \bar{z}}.
\]

Applying (\ref{e:geodesic}), the oriented geodesic $\gamma=(\xi,\eta)$ is
\[
\xi=\xi'\sinh(r_0+i\theta)\qquad\qquad \eta=\eta'+\frac{1}{\bar{\xi}'\tanh(r_0-i\theta)}
\]

Moving the point $p$ along the geodesic $(\xi',\eta')$, it is equivalent to an affine shift of $r_0$.
 
Therefore we obtain the surface $\Sigma$ given by the immersion $f:{\mathbb C}\rightarrow {\mathbb L}({\mathbb H}^3):(\nu,\bar{\nu})\mapsto (\xi(\nu,\bar{\nu}),\eta(\nu,\bar{\nu}))$ where
\[
\xi=\xi'\sinh\nu\qquad\qquad \eta=\eta'+\frac{1}{\bar{\xi}'\tanh\bar{\nu}},
\]
with $\nu=r+i\theta$. 

If we change the coordinates from $(\xi,\eta)$ to holomorphic coordinates $(\mu_1,\mu_2)$ on ${\mathbb{L}}({\mathbb{H}}^3)$, 
the surface $\Sigma$ is given by the following immersion
\[
\mu_1(\nu,\bar{\nu})=\frac{1-\cosh\bar{\nu}-\eta'\bar{\xi}'\sinh\bar{\nu}}{\bar{\xi}'\sinh\bar{\nu}}\qquad\qquad \mu_2(\nu,\bar{\nu})=\frac{\xi'\sinh\nu}{1+\cosh\nu+\bar{\eta}'\xi'\sinh\nu}.
\]

We can easily see that
\[
\sinh\nu=\frac{2(\xi')^{-1}\mu_2}{1+\bar{\mu}_1\mu_2}\qquad\qquad\mbox{and}\qquad\qquad \cosh\nu=\frac{1-\bar{\mu}_1\mu_2-2\bar{\eta}'\mu_2}{1+\bar{\mu}_1\mu_2},
\]
and from the identity $\cosh^2\nu-\sinh^2\nu=1$ we find that the Lagrangian surface $\Sigma$ is a maximal surface since it can be written 
\[
\bar{\mu}_2=\frac{\lambda_1\mu_1+1}{\mu_1+\lambda_2},\qquad\qquad\mbox{with}\qquad\qquad\lambda_1=\frac{1}{\eta'}\qquad \lambda_2=\frac{1}{\eta'}\left[(\eta')^2-\frac{1}{(\bar{\xi}')^2}\right],
\]
which completes the proof. 
\end{pf}

\vspace{0.1in}

\noindent
{\bf Note:} We have proved on the main theorem that every maximal Lagrangian surface is given by the graph (\ref{e:maximalimmersion}) and is orthogonal to a family of parallel equidistant tubes $\{S_{t}\}_{t\in I}$ around to the following oriented geodesics $\gamma'=(\mu'_1,\mu'_2)$ and $\tilde{\gamma}'=(\tilde{\mu}'_1,\tilde{\mu}'_2)$, given by
\[
\mu'_1=\frac{-1+\sqrt{1-\lambda_1\lambda_2}}{\lambda_1},\qquad \mu'_2=\frac{\bar{\lambda}_1}{1+\sqrt{1-\bar{\lambda}_1\bar{\lambda}_2}},
\]
and 
\[
\tilde{\mu}'_1=-\frac{1+\sqrt{1-\lambda_1\lambda_2}}{\lambda_1}\qquad \tilde{\mu}'_2=\frac{\bar{\lambda}_1}{1-\sqrt{1-\bar{\lambda}_1\bar{\lambda}_2}}.
\]
Consider now the antipodal map $\tau:{\mathbb P}^1\rightarrow {\mathbb P}^1:x\mapsto -\bar{x}^{-1}$ and observe that $\tilde{\mu}'_1=\tau(\mu'_2)$ and $\tilde{\mu}'_2=\tau(\mu'_1)$ which means that the geodesic $\tilde{\gamma}$ is obtained by reversing the orientation of the geodesic $\gamma$. In other words $\tilde{\gamma}$ and $\gamma$ describe the same geodesic, up to orientation.

\end{document}